\documentclass{article}
\usepackage[paperwidth=7in, paperheight=9.5in, margin=.875in]{geometry}

\usepackage{array}
\usepackage{cite}

\usepackage[caption=false]{subfig}
\usepackage{pgfplots}

\usepackage{algorithmic}

\usepackage{graphicx,epstopdf}

\usepackage{subfig}
\usepackage{amsmath}
\usepackage{amssymb}
\usepackage{graphicx}
\usepackage{dcolumn}
\usepackage{mathtools}
\usepackage{amsmath,amsfonts}
\usepackage{bm}
\usepackage{amsmath}
\usepackage{amssymb}
\usepackage{color}
\usepackage{float}
\usepackage{setspace}
\usepackage{tabularx}

\allowdisplaybreaks

\newcommand{\f}{\frac}
\newcommand{\ds}{\displaystyle}
\newcommand{\E}{ {\mathbb{E}} }

\newcommand{\be}{\begin{equation}}
\newcommand{\ee}{\end{equation}}

\colorlet{texcscolor}{blue!50!black}
\colorlet{texemcolor}{red!70!black}
\colorlet{texpreamble}{red!70!black}
\colorlet{codebackground}{black!25!white!25}

\date{}

\title{A Drift Homotopy Implicit Particle Filter Method for Nonlinear Filtering problems}

\author{
Xin Li \thanks{Department of Mathematics, Florida State University, Tallahassee, Florida.}
\and Feng Bao\thanks{ Department of Mathematics, Florida State University, Tallahassee, Florida, \ ({\tt fbao@fsu.edu}).}
\and Kyle Gallivan \thanks{Department of Mathematics, Florida State University, Tallahassee, Florida,  \ ({\tt kgallivan@fsu.edu}).}
      }

\begin{document}
\maketitle

\begin{abstract}
In this paper, we develop a drift homotopy implicit particle filter method. The methodology of our approach is to adopt the concept of drift homotopy in the resampling procedure of the particle filter method for solving the nonlinear filtering problem, and we introduce an implicit particle filter method to improve the efficiency of the drift homotopy resampling procedure.  Numerical experiments are carried out to demonstrate the effectiveness and efficiency of our drift homotopy implicit particle filter.
\end{abstract}
\vspace{0.5em}

\noindent \textit{\textbf{Keywords:} } \
\noindent Nonlinear filtering problem, particle filter,  drift homotopy dynamics, implicit sampling\\


\section{Introduction}

The goal of nonlinear filtering problems is to make the best estimation for the state of some stochastic dynamical system based on its partial noisy observations. As a key mathematical tool for data assimilation, nonlinear filtering methods have various applications in many scientific and engineering areas, such as weather forecasting, parameter estimation, signal processing, target tracking, and machine learning \cite{DA-applications}. 
\vspace{0.4em}

There are two types of approaches to solve the nonlinear filtering problem. The first type formulates the conditional distribution of the state of the target dynamical system by stochastic partial (or ordinary) differential equations, and then computes approximated distributions through numerical solutions of the equations \cite{zakai, time-discretization, Crisan-Zakai, Gobet-Zakai, Bao_Zakai, BSDE_filter, Bao_filter_jump, Bao_CMS, BCZ_2018, Bao_first, Bao_CiCP_2019, Bao_Semi, Bao_Half}. 
The other type of approach is called the ``Bayesian filter'', in which Bayesian inference is used to incorporate observational data into the dynamical model to estimate the state.  In this work, we focus on the Bayesian approach due to its wide acceptance by practitioners. There are two categories of Bayesian filter: Kalman-type filters and particle filters. The main idea of Kalman-type filers is to use the classic Kalman-Bucy filter, which can solve the linear filtering problem analytically, to solve a linearized filtering problem. Well-known  Kalman-type filters include the ensemble Kalman filter, the extended Kalman filter, and the unscented Kalman filter \cite{ExKF, UnKF, EvensenBook, Evense_EnKF, Tong_EnKF}. Although Kalman filters are easy to implement and they are efficient in solving high dimensional problems, the major drawback of the Kalman type filters is their poor performance when the dynamical system and the observation function are highly nonlinear. The linearized problem does not provide a good approximation to the original nonlinear filtering problem, and the quality of the state estimate degrades significantly. The particle filter method (also called the sequential Monte Carlo method), on the other hand, is designed specifically to solve nonlinear filtering problems. In contrast to the Kalman filter framework, which propagates and updates Gaussian distributions, the particle filter uses a set of weighted random samples (particles) to describe arbitrary conditional distributions of the target state \cite{cd2002, Doucet_PF, particle-filter, Bao_parameter, Kunsch_PF, BCH}. The weights on the particles are used to incorporate observational information, and they are assigned by Bayesian inference.  Through flexible Monte Carlo sampling, a particle filter can effectively handle nonlinear dynamics and nonlinear observations. The primary challenge of particle filter is that the errors caused by Monte Carlo sampling can accumulate through the sequential sampling procedures. Therefore, particle filters often suffer from the so-called ``degeneracy'' problem. In other words, after several estimation steps, most particles tend to lie in insignificant regions of the distribution, hence the effective ensemble size is dramatically reduced \cite{Sny-particle}. 

\vspace{0.4em}

One of the most effective approaches to address the degeneracy problem in the particle filter is ``resampling''. The goal of resampling is to rejuvenate the particle cloud and relocate particles from low density regions to high density regions.
Usually, a resampling step is implemented after (or combined with) Bayesian inference and it generates a set of particles (or moves current particles) that follow the desired conditional distribution of the target state. Successful resampling methods include sequential importance sampling with resampling (the benchmark method), the auxiliary particle filters, the Markov Chain Monte Carlo particle filter, the drift homotopy particle filter, and the implicit particle filter \cite{particle-filter, APF, Andrieu_MC, CT1, MCMC-PF, vanLeeuwen, Maroulas_PF}. While all of these methods can mitigate the degeneracy problem to some extent, they all have their disadvantages and drawbacks. 
\vspace{0.4em}

In this paper, we develop a drift homotopy implicit particle filter method that combines the drift homotopy particle filter \cite{Maroulas_PF, Kang-PF} and the implicit particle filter \cite{CT1}. The central concept of the drift homotopy particle filter is to construct a sequence of intermediate systems called drift homotopy dynamics, and then transport particles by using the Markov Chain Monte Carlo (MCMC) sampling method, which is driven by those intermediate homotopy systems, to high density regions of the desired state distribution.
The drift homotopy dynamics are usually designed in a way so that the observational data play a more important role in the first few  particle transportation steps, and then the original filtering dynamical model is incorporated into the state distribution gradually. In this way, the drift homotopy particle filter is different from most traditional Bayesian approaches, which simulate dynamical models first and then incorporate data through Bayesian inference. As a result, the drift homotopy systems provide a mechanism to ``trust the observational data'' first, which typically results more robust estimation performance. The main drawback of the drift homotopy particle filter is that the MCMC sampling procedure is time consuming -- especially carried out repeatedly through the sequence of drift homotopy dynamics , and thus the drift homotopy particle filter is not an efficient method. The primary effort of the implicit particle filter is to carry out an implicit sampling procedure, which works by first picking target probabilities and then looking for particles that assume them, so that the particles are guided to the high probability region. In practical implementation, the implicit sampling procedure is achieved by optimization, and the efficiency of implicit sampling is based on the shape of state distribution, which is governed by the state dynamics and the observational data. Therefore, when the state dynamical model and the observational data do not align well, the optimization task in the implicit particle filter could be very challenging. 
\vspace{0em}

Our motivation for developing a drift homotopy implicit particle filter (DHIPF) method is to exploit the advantages of the drift homotopy particle filter and the implicit particle filter while alleviating their disadvantages. Specifically, we shall adopt the general drift homotopy framework and utilize a sequence of drift homotopy dynamics to transport particles. However, instead of using the MCMC sampling method to move particles slowly, we treat the sampling procedure for intermediate drift homotopy systems as a sequence of nonlinear filtering problems and then use the implicit particle filter to solve those filtering problems. Since the implicit sampling in the implicit particle filter is achieved by optimization, the implicit particle filter sampler for intermediate drift homotopy systems is much more efficient than the MCMC sampling method. In this connection, the application of the implicit particle filter in DHIPF can significantly improve efficiency of the conventional drift homotopy particle filter. On the other hand, since the observational data play a more important role in the first few particle transportation steps in the drift homotopy systems, our DHIPF method could endow the implicit particle filter with the mechanism of ``trust observational data first'', which can make the DHIPF obtain more robust estimation results.

\vspace{0.4em}

The rest of this paper is organized as follows. In Section 2, we introduce the general mathematical formulation of the nonlinear filtering problem. In Section 3, we introduce the state-of-the-art approach to solve the filtering problem, i.e. the particle filter, and briefly discuss the drift homotopy particle filter and the implicit particle filter. Then, in Section 4, we combine the drift homotopy particle filter and the implicit particle filter to establish our drift homotopy implicit particle filter method. Numerical experiments that illustrate the performance of our method are given in Section 5.


\section{The nonlinear filtering problem}

An optimal filtering problem is usually described by a system of stochastic differential equations (SDEs)
\begin{equation}\label{NLF-continuous}
\begin{aligned}
d X_t =& f(X_t) dt + \sigma_t dW_t,  \vspace{4em}  &\text{(State)}\\
dY_t =& g(X_t) dt + d V_t.  &\text{(Observation)}
\end{aligned}
\end{equation}
The first equation in \eqref{NLF-continuous} is a state equation that models the state of a dynamical system, where $W_t$ is a standard Brownian motion, {\color{blue}$f(X_t) dt$ is the drift term} and the $\sigma_t dW_t$ integral is called the diffusion term. The second equation is the observational equation that gives partial noisy observations of $X_t$. In practical applications, this continuous version of the optimal filtering problem is often discretized and represented by the following discrete state-space model 
\begin{equation}\label{NLF}
\begin{aligned}
X_{n+1} =& f(X_n) + \sigma_n w_n,  \vspace{4em}  &\text{(State)}\\
Y_{n+1} =& g(X_{n+1}) + v_n,  &\text{(Observation)}
\end{aligned}
\end{equation}
where we have incorporated the temporal discretization steps $\Delta t$ into the model. In this way,
the sequence $\{X_n\}_n \in \mathbb{R}^d$ describes the state of the stochastic dynamical system, the function $f: \mathbb{R}^d \rightarrow \mathbb{R}^d$ now plays the role of the drift term in the continuous state equation in \eqref{NLF-continuous}, and the state of $X$ is perturbed by a sequence of $d$-dimensional standard Gaussian noises $\{w_n\}_n$ with their coefficients $\{\sigma_n\}_n$, and $Y_{n+1} \in \mathbb{R}^m$ is the $m$-dimensional partial noisy measurement on $X_{n+1}$ through the observation function $g: \mathbb{R}^d \rightarrow \mathbb{R}^m$, which is also perturbed by a Gaussian noise $v_n$ independent from $w_n$ with the standard deviation $R$. For the discretized optimal filtering problem \eqref{NLF}, the first stochastic process in \eqref{NLF} is called the ``state process'' and the second process is called the ``observation process''. When the functions $f$ and $g$ are nonlinear functions, the filtering problem is called the ``nonlinear filtering problem''. The goal of the nonlinear filtering problem is to find the best estimate for $\Phi(X_{n+1})$ given the observational data $Y_{1:n+1} := \{ Y_1, Y_2, \cdots, Y_{n+1}\}$, where $\Phi$ is a test function that represents the quantity of interest in the nonlinear filtering problem. Mathematically, we aim to find the ``optimal filter'' $\tilde{\Phi}(X_{n+1})$ as a conditional expectation given $Y_{1:n+1}$, i.e. 
$$\tilde{\Phi}(X_{n+1}): = \E[\Phi(X_{n+1}) | Y_{1:n+1}].$$

The standard approach to solve the optimal filtering problem \eqref{NLF} is the ``Bayesian filter'', which aims to find the best estimate for the conditional probability density function (pdf) $p(X_{n+1} | Y_{1:n+1})$ of the state through recursive Bayesian inference. Then, the conditional pdf $p(X_{n+1} | Y_{1:n+1})$, which is also called the ``filtering density'', can be used to calculate the optimal filter $\tilde{\Phi}$. Specifically, the Bayesian filter is composed of two steps: a prediction step and an update step. For the given conditional pdf $p(X_n|Y_{1:n})$ at the time instant $n$, the prediction step is carried out by the following Chapman-Kolmogorov formula,
\begin{equation}\label{C-K}
p(X_{n+1} | Y_{1:n}) = \int p_f(X_{n+1} | X_n)  p(X_n | Y_{1:n}) dX_n, \qquad \text{(Prediction)}
\end{equation} 
where $p_f(X_{n+1}|X_n)$ is the transition probability associated with the state dynamical function $f$, and $p(X_{n+1} | Y_{1:n})$ is called the prior pdf, which predicts the state of $X_{n+1}$ before receiving the new observational data. In the update step, we apply the Bayes formula to incorporate the observational data $Y_{n+1}$ to update the prior pdf and get a posterior distribution $p(X_{n+1}|Y_{1:n+1})$ as following
\begin{equation}\label{Bayes}
p(X_{n+1} | Y_{1:n+1}) = \f{p(Y_{n+1}| X_{n+1}) p(X_{n+1} | Y_{1:n})}{p(Y_{n+1} | Y_{1 : n})}, \qquad \text{(Update)}
\end{equation} 
where $p(Y_{n+1}| X_{n+1})$ is the likelihood function that measures the discrepancy between the state and the observation, and the denominator in \eqref{Bayes} is given by
$$p(Y_{n+1} | Y_{1 : n}) = \int p(Y_{n+1} | X_{n+1}) p(X_{n+1} | Y_{1:n}) d X_{n+1},$$
which normalizes the posterior $p(X_{n+1} | Y_{1:n+1})$. 

In this work, we aim to develop an efficient and effective particle filter method to implement the Bayesian filter's ``Prediction-Update'' framework.

\section{The particle filter approach}

\subsection{The generic particle filter framework}

The main strategy of the particle filter is to use a cloud of samples, which are called particles, to represent conditional distributions, and use the recursive Bayesian filter framework \eqref{C-K} - \eqref{Bayes} to propagate and updates the particle cloud.
In what follows, we give a brief discussion to introduce the general framework of particle filters.

At the time instant $n$, assume that we have a set of $N_p$ particles that form an empirical distribution $\tilde{p}(X_{n} | Y_{1:n})$, which approximates the conditional pdf $p(X_{n} | Y_{1:n})$. We denote this set of particles by $\{x_n^{(i)}\}_{i=1}^{N_p}$, and the empirical distribution is defined as
\begin{equation}\label{empirical}
\tilde{p}(X_{n} | Y_{1:n}) : = \f{1}{N_p} \sum_{i=1}^{N_p} \delta_{x_n^{(i)}}(X_n),
\end{equation}
where $\delta_{x}$ is the Dirac delta function. In the prediction step, we propagate each sample $x_n^{(i)}$ in the particle cloud through the state dynamical model $f$ to get a predicted sample $\tilde{x}_{n+1}^{(i)}$. 
In this way, the ensemble of predicted particles $\{\tilde{x}_{n+1}^{(i)}\}_{i=1}^{N_p}$ form an empirical distribution $\tilde{\pi}(X_{n+1} | Y_{1:n})$ defined as
\begin{equation}\label{empirical-pi}
\tilde{\pi}(X_{n+1} | Y_{1:n}) : = \f{1}{N_p} \sum_{i=1}^{N_p}\delta_{\tilde{x}_{n+1}^{(i)}}(X_{n+1}),
\end{equation}
which is an approximation for the prior distribution $p(X_{n+1} | Y_{1:n})$.

In the update step, after receiving the new measurement $Y_{n+1}$, we carry out Bayesian inference through the Bayes formula \eqref{Bayes} to incorporate the new observational data to get an approximation for the posterior distribution. Specifically, we use the empirical distribution $\tilde{\pi}(X_{n+1}|Y_{1:n})$ to replace the prior distribution $p(X_{n+1}|Y_{1:n})$ in \eqref{Bayes}, and obtain
\begin{equation}\label{update:integral}
\tilde{\pi}(X_{n+1} | Y_{1:n+1}) = \f{p(Y_{n+1} | X_{n+1}) \tilde{\pi}(X_{n+1}|Y_{1:n})}{\int p(Y_{n+1 | X_{n+1}})  \tilde{\pi}(X_{n+1}|Y_{1:n}) d X_{n+1}},
\end{equation}
where the likelihood function $p(Y_{n+1} | X_{n+1})$ of the Gaussian noise $v$ with the standard deviation $R$ is given by 
$p(Y_{n+1} | X_{n+1}) = \f{1}{\sqrt{(2 \pi R^2)^m}} \exp\Big( - \f{\| g(X_{n+1}) - Y_{n+1} \| ^2}{2 R^2}\Big)$.
With the empirical distribution $\tilde{\pi}(X_{n+1}|Y_{1:n})$ defined in \eqref{empirical-pi}, the Bayesian inference formula \eqref{update:integral} can be implemented by the following update scheme
\begin{equation}\label{PF-update}
\tilde{\pi}(X_{n+1} | Y_{1:n+1}) = \f{ \sum_{i=1}^{N_p} p(Y_{n+1} | \tilde{x}_{n+1}^{(i)}) \delta_{\tilde{x}_{n+1}^{(i)}}(X_{n+1})}{ \sum_{i=1}^{N_p} p(Y_{n+1} | \tilde{x}_{n+1}^{(i)}) }.
\end{equation}
We let $\alpha_{n+1}^{(i)} \propto p(Y_{n+1} | \tilde{x}_{n+1}^{(i)})$ be the importance weight corresponding to the particle $\tilde{x}_{n+1}^{(i)}$ such that $\ds \sum_{i=1}^{N_p} \alpha_{n+1}^{(i)} = 1$.
Then, the empirical distribution $\tilde{\pi}(X_{n+1} | Y_{1:n+1})$ is an approximation for the posterior distribution $p(X_{n+1} | Y_{1:n+1})$, where we have 
$$\tilde{\pi}(X_{n+1} | Y_{1:n+1}) = \sum_{i=1}^{N_p} \alpha_{n+1}^{(i)} \delta_{\tilde{x}_{n+1}^{(i)}}(X_{n+1}).$$

In practice, due to the extra uncertainties involved in the observational data and the sequential sampling errors, the weights on many particles tend to be negligible after several recursive steps, and only a few particles have very large weights, which significantly reduces the effective ensemble size. This fact of losing the effectiveness of particle sizes is often called the ``degeneracy'' of particles. To address the degeneracy problem, a \textit{resampling step} is introduced to re-generate the particle cloud with equally weighted particles that describe the empirical distribution $\tilde{\pi}(X_{n+1} | Y_{1:n+1})$. In the benchmark bootstrap particle filter, which is also known as the sequential importance sampling with resampling method, people use importance sampling to generate $N_p$ samples, denoted by $\{x_{n+1}^{(i)}\}_{i=1}^{N_p}$, which include more copies of particles in the weighted particle cloud $\{\tilde{x}_{n+1}^{(i)}\}_{i=1}^{N_p}$ and discard the low weight ones. In this way, the resampled particles give us the following empirical distribution
\begin{equation*}
\tilde{p}(X_{n+1} | Y_{1:n+1}) : = \f{1}{N_p} \sum_{i=1}^{N_p} \delta_{x_{n+1}^{(i)}}(X_{n+1}).
\end{equation*}

Although the importance sampling in the bootstrap particle filter is a successful strategy to increase the effective particle size in recursion, simply reproducing more particles in the original predicted particle set $\{\tilde{x}_{n+1}^{(i)}\}_{i=1}^{N_p}$ still suffers from the degeneracy problem since the prediction step is not informed by the observational information and the predicted particles may not provide good candidate particles to be reproduced. To conquer the degeneracy of particles, several data informed resampling methods have been developed, such as the auxiliary transportation, the Markov Chain Monte Carlo resampling, the drift homotopy particle filter, the implicit particle filter. In this work, we combine the drift homotopy particle filter and the implicit particle filter to develop an efficient drift homotopy implicit particle filter method that takes the advantages of both methods. To proceed, we shall introduce the drift homotopy particle filter and the implicit particle filter in the following.

\subsection{The drift homotopy particle filter}

Homotopy is a mathematical process that continuously transforms one function to another. When adopting homotopy in the particle filter, we design a homotopy process that transforms the drift term in the state dynamics gradually from an intermediate drift homotopy system, so that the sampling of the desired filtering density can be implemented effectively. 

In this connection, the key of the drift homotopy particle filter is to construct a sequence of stochastic dynamical systems with modified drift terms that interpolate between the original and modified drifts in the state model.
To proceed, assume that we have an empirical distribution that describes the conditional pdf $p(X_n | Y_{1:n})$  at the time instant $n$ with equally weighted particles $\{x_{n}^{(i)}\}_{i=1}^{N_p}$. For a given particle $x_{n}^{(i)}$, the ``ideal'' particle at the time instant $n+1$ that could be used to describe the posterior distribution $p(X_{n+1} | Y_{1:n+1})$ should follow the distribution $p(x_{n+1} | x_{n}^{(i)}, Y_{1:n+1})$, which is defined by
\begin{equation}\label{Posterior-original}
p(X_{n+1} | x_{n}^{(i)}, Y_{1:n+1}) \propto p_f(X_{n+1} | x_{n}^{(i)}) p(Y_{n+1}|x_{n+1}),
\end{equation}
where $p_f(X_{n+1} | x_{n}^{(i)})$ is the transition probability that describes the conditional distribution of $X_{n+1}$ with the given initial state at the time instant $n$ as $X_n = x_n^{(i)}$, and the likelihood $p(Y_{n+1} | X_{n+1})$ incorporates the observational data $Y_{n+1}$ into the conditional distribution $p(X_{n+1} | x_{n}^{(i)}, Y_{1:n+1})$. In order to generate a sample that follows the posterior distribution, which considers the new observational data $Y_{1:n+1}$, instead of using the importance sampling method which only reproduces more copies of high density propagation samples, one could use the Markov Chain Monte Carlo (MCMC) sampling method to generate a sample $x_{n+1}^{(i)}$ that follows the conditional distribution $p(X_{n+1} | x_{n}^{(i)}, Y_{1:n+1})$. However, it is well known that the effectiveness and efficiency of MCMC sampling depends on the complexity of the target distribution. When the dynamics of the state process is complicated, it is difficult for the MCMC method to generate the sample $x_{n+1}^{(i)}$. Moreover, in the case that the practical observational data has large deviation from the prediction, the resulting posterior distribution might have large covariance, which makes the MCMC procedure have even worse performance.   

The primary contribution of the drift homotopy particle filter is to improve the performance of the MCMC procedure to effectively generate a sample that follows the posterior distribution. Also, it could provide a mechanism that allows us to ``trust the data'' first. In the drift homotopy particle filter, instead of generating a sample from $p(X_{n+1} | x_{n}^{(i)}, Y_{1:n+1})$ by using MCMC sampling directly, we introduce a sequence of drift homotopy dynamical systems given as follows
\begin{equation}\label{DH-state}
X^{\ast} = (1 - \beta_l) b(X) + \beta_l f(X) + \sigma_n w_n, \qquad l = 0, 1, 2, \cdots, L
\end{equation} 
where the function $b$ is artificially defined intermediate drift term, which is different from the original drift $f$, and $\{\beta_l\}_{l=0}^{L}$ is a constant sequence increasing from $0$ to $1$. Therefore, when $l=0$, we have $\beta_0 = 0$, and the dynamical system \eqref{DH-state} only contains the intermediate function $b$. On the other hand, when $l = L$, i.e. $\beta_L = 1$, the intermediate drift term $b$ is gone and \eqref{DH-state} coincides the original state dynamical model in the nonlinear filtering problem \eqref{NLF}. For a specific drift homotopy level $l$, $ 0 \leq l \leq L - 1$, the dynamical system \eqref{DH-state} is driven by the combined model $(1 - \beta_l) b(X^{l}) + \beta_l f(X^{l})$. To generate a sample that follows the desired conditional distribution $p(X_{n+1} | x_{n}^{(i)}, Y_{1:n+1})$ (as described in \eqref{Posterior-original}), we incorporate the observational data $Y_{n+1}$ through the likelihood function of the dynamical system \eqref{DH-state}.  Since the target distribution $p(X_{n+1} | x_{n}^{(i)}, Y_{1:n+1})$ is conditioned on $x_n^{(i)}$, we take the state $X$ of the drift homotopy dynamics \eqref{DH-state} to be the particle $x_n^{(i)}$, i.e. $X = x_n^{(i)}$, which is a sample in the previous conditional particle cloud $\{x_n^{(i)}\}_{i=1}^{N_p}$. As a result, we obtain a sequence of drift homotopy distributions $\{p_l(X^{\ast} | x_n^{(i)}, Y_{1:n+1})\}_{l=1}^{L}$ for each particle $x_n^{(i)}$ through the following proportional relation
\begin{equation}\label{Posterior:DH}
p_l(X^{\ast}| x_n^{(i)}, Y_{1:n+1}) \propto p_l(X^{\ast} | x_n^{(i)}) p(Y_{n+1}|X^{\ast}),
\end{equation}
where $p_l(X^{\ast} | x_n^{(i)}) $ is the transition probability of the system \eqref{DH-state} determined by the combined dynamics. In this way, the distribution of the $L$-th drift homotopy step is the desired posterior distribution given that $p_L(X^{\ast} | x_n^{(i)}) = p_f(X_{n+1} | x_n^{(i)})$.


To carry out the drift homotopy particle filter, we first choose $x_n^{(i)}$ as the initial state of $X$ and use the MCMC method to generate a sample, denoted by $x^{\ast, 0}_{n+1}$, from the first homotopy distribution $p_0(X^{\ast}| x_n^{(i)}, Y_{1:n+1})$. For the $l$-th drift homotopy level, $l = 0, 1, 2, \cdots, L-1$, assume that we have the sample $x^{\ast, l}_{n+1}$, we let $x^{\ast, l+1}_{n+1}$ be the initial state of the drift homotopy system \eqref{DH-state} and use the MCMC method to generate a sample, denoted by $x^{\ast, l+1}_{n+1}$ from the homotopy distribution $p_{l+1}(X^{\ast}| x_n^{(i)}, Y_{1:n+1})$. As a result, the sample $x^{\ast, L}_{n+1}$ that we generate in the iterative drift homotopy procedures gives a sample that follows the desired posterior distribution $p(X_{n+1} | x_{n}^{(i)}, Y_{1:n+1})$, and we let $x_{n+1}^{(i)} = x^{\ast, L}_{n+1}$. 

From the above discussion, we can see that the main theme of the drift homotopy particle filter is to transport the particle $x_n^{(i)}$ gradually to $x_{n+1}^{(i)}$, which is then used to formulate the posterior distribution. 
In order to fully assimilate the observational information into the posterior distribution, we choose the intermediate drift term $b$ so that the initial transition probability $p_0(X^{\ast}| x_n^{(i)})$ is well-aligned with the likelihood function $p(Y_{n+1}|X^{\ast})$.
Then, as the drift homotopy distributions $\{p_{l}(X^{\ast}| x_n^{(i)}, Y_{1:n+1})\}_l$ morph gradually to the ultimate posterior distribution $p(X_{n+1} | x_{n}^{(i)}, Y_{1:n+1})$, the original state dynamics $f$ is incorporated. In this way, the homotopy procedure builds a bridge to connect the sample $x_n^{(i)}$ in the previous particle cloud at the time instant $n$ to the particle $x_{n+1}^{(i)}$ for the desired posterior distribution. Moreover, note that the combined dynamical system only contains the intermediate drift term $b$ in the first drift homotopy step, which is designed to be well-aligned with the  observational data. Therefore, the observational data would have more influence in the drift homotopy distribution. In this way, the drift homotopy particle filter also creates a mechanism that allows us to ``trust the data first''.



\subsection{The implicit particle filter}

The central concept of the implicit particle filter is to update the particles by first deriving implicit transportation probabilities, which construct the conditional distribution of the target, and then determine particle locations that assume them. As a result, the implicit sampling procedure guides the particles one by one to the high probability domain of the desired posterior distribution, and therefore it can effectively address the degeneracy problem of the particle filter approach.  

In the standard particle filter, assume that we have the particle $x_{n}^{(i)}$ at the time step $n$, then the predicted particle location $\tilde{x}_{n+1}^{(i)}$ is propagated from the particle $x_{n}^{(i)}$ through the state dynamics introduced in \eqref{NLF}. According to the Bayesian update scheme \eqref{Bayes}, the weight on $\tilde{x}_{n+1}^{(i)}$ is given by
$$\alpha_{n+1}^{(i)} = \f{p(\tilde{x}_{n+1}^{(i)} | x_{n}^{(i)}) p(Y_{n+1} | \tilde{x}_{n+1}^{(i)})}{C},$$
where $C$ is the normalization constant. The set of weights $\{\alpha_{n+1}^{(i)}\}_{i=1}^{N_p}$, together with the predicted particle locations $\{\tilde{x}_{n+1}^{(i)} \}_{i=1}^{N_p}$, now describe the desired posterior density. In the implicit particle filter method, instead of simply reproducing more copies of particles in the predicted particle cloud with high weights, an implicit sampling procedure is used to generate ``optimal'' particle locations, i.e. $\{\hat{x}_{n+1}^{(i)} \}_{i=1}^{N_p}$, that follow the posterior distribution directly. 

To achieve this goal, we first choose a reference random variable $\xi$ with a pre-determined pdf that is easy to sample. Then, we treat the desired optimal sample variable $\hat{x}_{n+1}^{(i)} : = \psi^{n+1, i}(\xi)$ as a function of $\xi$, which is indexed by both the time instant and the specific choice of particle. The purpose of the function $\psi$ is to connect highly probable values of $\xi$ to highly probable values of $\hat{x}_{n+1}^{(i)}$, which follow the posterior distribution. To obtain the mapping $\psi$, we define a function $F^{n+1, i}$ corresponding to each time instant $n+1$ and each particle $i$, such that
\begin{equation}\label{Def:F}
\exp\Big(-F^{n+1, i} (\hat{x}_{n+1}^{(i)}) \Big) : = p(\hat{x}_{n+1}^{(i)} | x_{n}^{(i)}) p(Y_{n+1} | \hat{x}_{n+1}^{(i)}),
\end{equation}
and solve the equation
\begin{equation}\label{F:optimization}
F^{n+1, i}\Big(\psi^{n+1, i}(\xi) \Big) - \gamma^{n+1, i} = \f{\xi^{T} \xi}{2}
\end{equation}
to get the function $\psi^{n+1, i}$. The random variable $\gamma^{n+1, i}$ in the above equation is an additive factor introduced to make the above equation solvable, and it is typical choose $\gamma^{n+1, i} = \min F^{n+1, i} + \lambda$, where $\lambda$ represents a small perturbation. When implementing the implicit particle filter numerically, optimization type numerical solvers are needed to calculate $\psi^{n+1, i}$ through \eqref{F:optimization}.

Once the function $\psi^{n+1, i}$ is determined, we can obtain the particle set $\{ \hat{x}_{n+1}^{(i)}\}_{i=1}^{N_p}$ since the position of the particle $\hat{x}_{n+1}^{(i)}$ appears with the (unnormalized) probability $\ds \exp\big(-\f{\xi^{T} \xi}{2}\big) J_{\psi^{n+1, i}}^{-1}$, where $J_{\psi^{n+1, i}}$ denotes the Jacobian of $\psi^{n+1, i}$, and the weight on $\hat{x}_{n+1}^{(i)}$ equals $\ds \exp\big(-\f{\xi^{T} \xi}{2} \big) \exp(-\gamma^{n+1, i})$ \cite{MTAC2012}, i.e. 
\begin{equation}\label{IPF:weight}
\hat{\alpha}_{n+1}^{(i)} := \ds \exp\big(-\f{\xi^{T} \xi}{2} \big) \exp(-\gamma^{n+1, i}).
\end{equation}
 To generate equally weighted particles, we apply the importance sampling method to resample particles $\{\hat{x}_{n+1}^{(i)}\}_{i=1}^{N_p}$. Specifically, we normalize the weights $\{\hat{\alpha}_{n+1}^{(i)}\}_{i=1}^{N_p}$ obtained in \eqref{IPF:weight} so that $\sum_{i=1}^{N_p} \hat{\alpha}_{n+1}^{(i)} = 1$. For each of $N_p$ random numbers $\zeta_k$, $k = 1, 2, \cdots, N_p$, drawn from the uniform distribution on $[0, 1]$, we choose a point $x_{n+1}^{(k)}$ randomly from the particle set $\{\hat{x}_{n+1}^{(i)}\}_{i=1}^{N_p}$ such that 
 $$\sum_{j=1}^{k-1} \hat{\alpha}_{n+1}^{(j)} < \zeta_k < \sum_{j=1}^{k} \hat{\alpha}_{n+1}^{(j)}. $$
Then the particle set $\{ x_{n+1}^{(i)} \}_{i=1}^{N_p}$ follows the importance distribution described by the weighted sample pairs $\{ (\hat{x}_{n+1}^{(i)}, \hat{\alpha}_{n+1}^{(i)}) \}_{i=1}^{N_p}$, and each particle $x_{n+1}^{(i)}$ has an equal weight. Note that  the particle set $\{ x_{n+1}^{(i)} \}_{i=1}^{N_p}$ resampled from $\{ \hat{x}_{n+1}^{(i)} \}_{i=1}^{N_p}$ describes the conditional distribution of the target better than the particles resampled from $\{\tilde{x}_{n+1}^{(i)} \}_{i=1}^{N_p}$ since $\{ \hat{x}_{n+1}^{(i)} \}_{i=1}^{N_p}$ are already in the high density region through the implicit sampling.

Theoretically, the implicit particle filter method can generate equally weighted particles that follow the desired posterior distribution if the (nonlinear) equation \eqref{F:optimization} can be effectively solved, and the performance of the implicit particle filter is based on the performance of the optimization procedure that solves the equation. However, the complexity of the equation \eqref{F:optimization} depends on the function $F$, which relies on both the transition probability and the likelihood function (as indicated in \eqref{Def:F}). In this way, the implicit particle filter could be computationally expensive and challenging when the dynamical model is not well-aligned with the observational data.

\section{Drift homotopy implicit particle filter}

The main effort of this work is to combine advantages of the drift homotopy particle filter and the implicit particle filter to construct a drift homotopy implicit particle filter (DHIPF) method, which can effectively use observational data and efficiently generate particles that follow the filtering density of the target. 
The general framework of our approach adopts the drift homotopy procedure in the drift homotopy particle filter. Instead of using MCMC as a sampling method in the drift homotopy particle filter, we carry out the ``implicit sampling'' procedure introduced the implicit particle filter to generate the desired particles efficiently.  

To proceed, we recall that in the drift homotopy particle filter, the drift homotopy sequence \eqref{DH-state} builds a bridge that connects the intermedia dynamics $b$ to the original state dynamics $f$.   With the observational data incorporated through likelihood (as described in \eqref{Posterior:DH}), the drift homotopy distribution $p_l(X^{\ast}| x_n^{(i)}, Y_{1:n+1})$ in each drift homotopy step is proportional to $p_l(X^{\ast} | x_n^{(i)}) p(Y_{n+1}|X^{\ast})$.
In the DHIPF, instead of simply using the drift homotopy dynamics as a bridge to transport samples, we consider the random variable $X$ in the drift homotopy dynamics \eqref{DH-state} as the state of a nonlinear filtering problem at the time instant $n$ and consider $X^{\ast}$ as the state at the time instant $n+1$. Therefore, the drift homotopy distribution $p_l(X^{\ast}| x_n^{(i)}, Y_{1:n+1})$ is equivalent to the filtering density of the following nonlinear filtering problem
\begin{equation}\label{DH-filtering}
\begin{aligned}
\tilde{X}_{n+1} =& (1 - \beta_l) b(\tilde{X}_n) + \beta_l f(\tilde{X}_n) + \sigma_n w_n, &\text{(State)}\\
\tilde{Y}_{n+1} =& g(\tilde{X}_{n+1}) + v_n,  &\text{(Observation)}
\end{aligned}
\end{equation} 
given that the state $\tilde{X}_n$ is chosen as a particle $x_n^{(i)}$  in the previous particle cloud $\{x_n^{(i)}\}_{i=1}^{N_p}$ of the original nonlinear filtering problem \eqref{NLF} and the observation $\tilde{Y}_{n+1}$ is taken as the observational data $Y_{n+1}$. In other words, we have 
$$p_l(X^{\ast}| x_n^{(i)}, Y_{1:n+1}) = p_l(\tilde{X}_{n+1}| \tilde{X}_{n}, \tilde{Y}_{n+1})\big|_{\tilde{X}_{n} = x_n^{(i)}, \tilde{Y}_{n+1} = Y_{1:n+1}}.$$ 
Then, the implicit particle filter method can be applied to solve the nonlinear filtering problem \eqref{DH-filtering} and produce a particle that follows the drift homotopy distribution $p_l(X^{\ast}| x_n^{(i)}, Y_{1:n+1})$. 

Specifically, for an appropriately chosen reference random variable $\xi$, we solve for the function $\psi^{n+1, i}_l(\xi)$ that connects highly probable values of $\xi$ to highly probable values of $\hat{x}_{n+1, l}^{(i)}$, i.e. $\hat{x}_{n+1, l}^{(i)} = \psi^{n+1, i}_l(\xi)$, where $\hat{x}_{n+1, l}^{(i)}$ is a particle that follows the distribution $p_l(\tilde{X}_{n+1}| \tilde{X}_{n}, \tilde{Y}_{n+1})\big|_{\tilde{X}_{n} = x_n^{(i)}, \tilde{Y}_{n+1} = Y_{1:n+1}}$. To this end, we define a function $F^{n+1, i}_l$ (corresponding to the particle $x^{(i)}_n$) by
\begin{equation}\label{Def:F_l}
\exp\Big(-F^{n+1, i}_l (\hat{x}_{n+1, l}^{(i)}) \Big) : = p_l(\hat{x}_{n+1, l}^{(i)} | x_{n}^{(i)}) p(Y_{n+1} | \hat{x}_{n+1, l}^{(i)}),
\end{equation}
where $p_l$ is the transition probability of the $l$-th drift homotopy dynamics. Then, we solve the following equation
\begin{equation}\label{F_l:optimization}
F^{n+1, i}_l\Big(\psi^{n+1, i}_l(\xi) \Big) - \gamma_l^{n+1, i} = \f{\xi^{T} \xi}{2}
\end{equation}
to obtain the function $\psi^{n+1, i}_l$, where $\gamma_l^{n+1, i} := \min F_l^{n+1, i} + \lambda$ is the factor that makes \eqref{F_l:optimization} solvable as we introduced in the equation \eqref{F:optimization}. Then, with the connection function $\psi^{n+1, i}_l$ solved through the equation \eqref{F_l:optimization}, we can generate the position of the particle $\hat{x}_{n+1, l}^{(i)}$ through the expression $\ds \exp\big(-\f{\xi^{T} \xi}{2}\big) J_{\psi_l^{n+1, i}}^{-1}$ for a pre-chosen sample that follows $\xi$.

Similar to the procedure that solves the equation \eqref{F:optimization} in the implicit particle filter, we use an optimization-based approach to determine $\psi^{n+1, i}_l$ (for a given sample drawn from $\xi$). Note that the drift homotopy dynamics morph gradually from the intermediate drift $b$ to the original dynamical model $f$. Therefore, the transition probabilities $\{p_l\}_{l=1}^{L}$ between two successive drift homotopy levels have similar distributions. Hence the values of implicit functions $\{\psi^{n+1, i}_l\}_{l=1}^L$ should be close for two successive drift homotopy levels.
To take the advantage of those bridging drift homotopy dynamics, we use the sample $\hat{x}^{(i)}_{n+1,l-1}$ obtained in the $l-1$-th drift homotopy level as the initial condition for the optimization procedure when solving for $\psi^{n+1, i}_{l}$. As a result, the optimization for solving $\psi^{n+1, i}_{l}$ convergences quickly due to the ``good'' initial condition and the implicit particle filter can be carried out efficiently. 
 
In the last homotopy level $L$, the drift homotopy dynamics become the original state dynamics $f$. Therefore, once we obtain the function $\psi^{n+1, i}_L$ through the implicit sampling procedure \eqref{Def:F_l} - \eqref{F_l:optimization}, we obtain the sample $\hat{x}^{(i)}_{n+1, L} $ that follows the desired filtering density $p(X_{n+1} | x_{n}^{(i)}, Y_{1:n+1})$ given the particle $_n^{(i)}$ and the fact $p_L(X^{\ast} | x_n^{(i)}) = p_f(X_{n+1} | x_n^{(i)})$.
\vspace{1em}

Our DHIPF algorithm is summarized in Table \ref{Algorithm}.
\begin{table}[h!]\caption{}\label{Algorithm}
\vspace{1em}
\centering
\begin{tabular} {p{0.9\textwidth}}
\hline\noalign{\smallskip}
{\bf Algorithm}: {\em Drift homotopy implicit particle filter (DHIPF) }\\
\noalign
{\smallskip}\hline
\noalign{\smallskip}
\vspace{-0.1cm}
\begin{spacing}{1.1}
\begin{algorithmic}\label{algorithm}
\item[Initialize] the particle cloud $\{x_0^{(i)}\}_{i=1}^{N_p}$, the number of drift homotopy levels $L$ with the intermediate drift function $b$ and the constant sequence $\{\beta_l\}_{l=0}^{L}$, and the reference random variable $\xi$ for the implicit particle filter procedure. \\
\vspace{0.25em}
\item[\textbf{while}] $n =0, 1, 2, \cdots $, \textbf{do} 
\vspace{-0.75em}
\begin{description}
\item[\textbf{for}] particles $i = 1, 2, \cdots, N_p$,   
\vspace{-1em}
\begin{description}
\item[\textbf{for}] drift homotopy levels $l = 0, 1, 2, \cdots, L-1$,   
\vspace{-1em}
\begin{description}
\item[-]  Construct the drift homotopy dynamics \eqref{DH-state};\vspace{-0.75em}
\item[-]  Solve for $\psi_{l}^{n+1, i}$ in the equation \eqref{F_l:optimization} with the initial guess $\hat{x}_{n+1, l}^{(i)}$;
\vspace{-1.25em}
\item[-] Generate the sample $\hat{x}_{n+1, l+1}^{(i)}$ through $\ds \exp\big(-\f{\xi^{T} \xi}{2}\big) J_{\psi_l^{n+1, i}}^{-1}$;
\vspace{-1em}
\end{description}
\item[\textbf{end for}]
\end{description}
\item[\textbf{end for}]
\vspace{-.75em}
\end{description}
\vspace{-.5em}
\item[] The particles $\{x_{n+1}^{(i)}\}_{i=1}^{N_p} := \{\hat{x}_{n+1, L}^{(i)}\}_{i=1}^{N_p}$ provide an empirical distribution for the filtering density  $p(X_{n+1} | Y_{1:n+1})$
\vspace{0.25em}
\item[\textbf{end while}]
\end{algorithmic}
\vspace{-1.2em}
\end{spacing}\\
\hline
\end{tabular}
\end{table}
Based on the above discussions for DHIPF, we can see that the optimization based implicit sampling procedure can generate the drift homotopy sample $\hat{x}_l^{(i)}$ much more efficiently -- compared with the standard MCMC sampling method.
On the other hand, the drift homotopy procedure creates a bridge that connects an intermediate dynamical function and the original state dynamics. Since two successive drift homotopy dynamics are similar, solutions of drift homotopy filtering problems \eqref{DH-filtering} change gradually to the filtering density of the original nonlinear filtering problem. Therefore, it is easy to achieve the optimality condition in the implicit particle filter , and hence the implicit particle filter can be implemented efficiently under our DHIPF framework. Moreover, the intermediate dynamics $b$ in the DHIPF is designed in a way so that the likelihood function would dominate the first few drift homotopy steps. Then, the drift homotopy procedure incorporates the original filtering dynamics and let the dynamical model combine with the observational data. In other words, the drift homotopy procedure aims to construct the desired filtering density starting from the likelihood instead of starting from the predicated model, which is typically implemented by most optimal filtering methods.  In this way, the DHIPF could endow the implicit particle filter the mechanism that trusts the observational data first.


\section{Numerical experiments}
In this section, we present two benchmark numerical examples to demonstrate the performance of our DHIPF method. In the first example, we track the state of a stochastic dynamical system driven by the double well potential. To demonstrate the advantageous performance of DHIPF, we compare our method with the implicit particle filter and the standard drift homotopy particle filter (with MCMC sampling ) -- along with other ``state-of-the-art'' methods. In the second example, we solve a Lorenz attractor problem. The target state that we estimate is driven by the Lorenz 63 model, which is a well-known chaotic dynamical system. We show that our method can capture the unpredicted chaotic behavior of the model by effective processing of observational data. All the numerical experiments are carried out on an Intel Core i7-670HQ 2.6GHz CPU.

\subsection{Double well potential}
The double well potential is an important quartic model in quantum mechanics and quantum field theory, and models derived from the double well potential have been widely used in nano-phase materials \cite{Bao_Acta}. The potential $U$, described by
\begin{equation*}
U(x) = \f{\alpha}{4} (x^4 - 2 x^2), 
\end{equation*}
has two stable positions at $x = 1$ and $x = -1$, where $\alpha$ is the model parameter that determines the ``depth'' of potential wells. While a particle is at positions other than $1$ and $-1$, it will be pushed by a force with the magnitude of $U'(x)$ towards one of the stable positions. In this example, we aim to estimate the state of a stochastic dynamical system driven by the double well potential, i.e.
\begin{equation*}
d X_{t} =- \alpha (X_{t}^3 - X_t) dt + \sigma dW_t, \hspace{2em} 0 \leq t \leq T 
\end{equation*}
and the data that we use to track the target state are direct observations on $X$, which are perturbed by Gaussian noises with standard deviation $R$.

We consider the following discretized nonlinear filtering problem
\begin{equation}\label{Ex1:NLF}
\begin{aligned}
X_{n+1} =& X_n - \alpha (X_{n}^3 - X_n) \Delta t + \sigma w_n,  \\
Y_{n+1} =& X_{n+1} + R v_n, 
\end{aligned}
\end{equation}
where $w_n$ and $v_n$ are two independent Gaussian random variables, and we track the state of $X$ for $300$ time steps with stepsize $\Delta t = 0.01$. In this example, we compare our DHIPF with four most successful nonlinear filtering methods: the auxiliary particle filter (APF), the ensemble Kalman filter (EnKF), the implicit particle filter and the drift homotopy particle filter (DHPF), where IPF and DHPF (implemented by MCMC sampling) are fundamental components that we use to construct our DHIPF.  For all the particle filters, we use $20$ particles to describe the one-dimensional state distribution and we use an ensemble of {\color{blue}$200$ Kalman filter samples} in the EnKF. Also, we choose the total number of drift homotopy steps to be $L = 2$, i.e. we use three intermediate dynamical systems to transport particles.

To provide a comprehensive demonstration of the performance of all the nonlinear filtering methods, we solve the double well potential tracking problem for three different cases.

\subsubsection*{Case 1.}
In this case, we choose the parameters for the double well potential tracking problem as $\alpha = 1$, $\sigma = 1.5$, and $R = 1.5$, and the initial state $X$ is set to be $X_0 = 0.6$.
\begin{figure}[h!]
\begin{center}
\subfloat[Tracking performance. ]{\includegraphics[scale = 0.55]{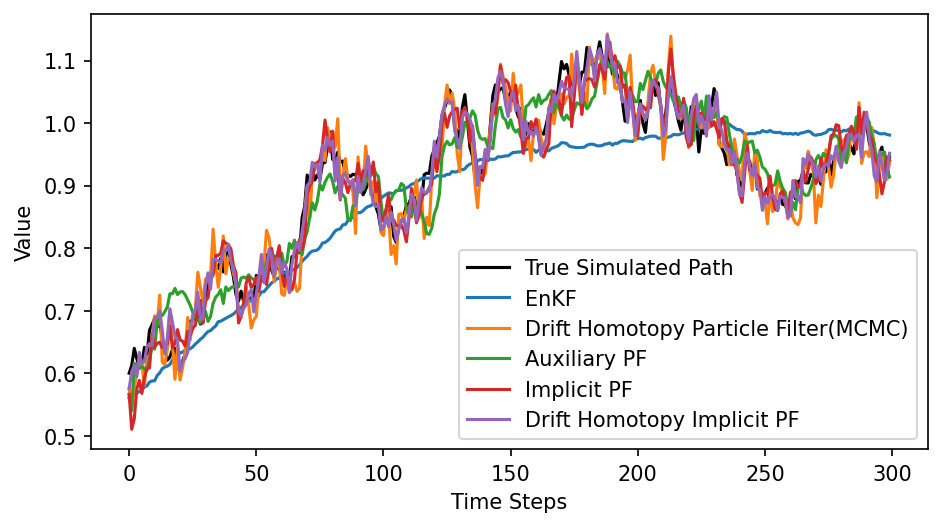} } \\ 
\subfloat[Estimation errors.]{\includegraphics[scale = 0.55]{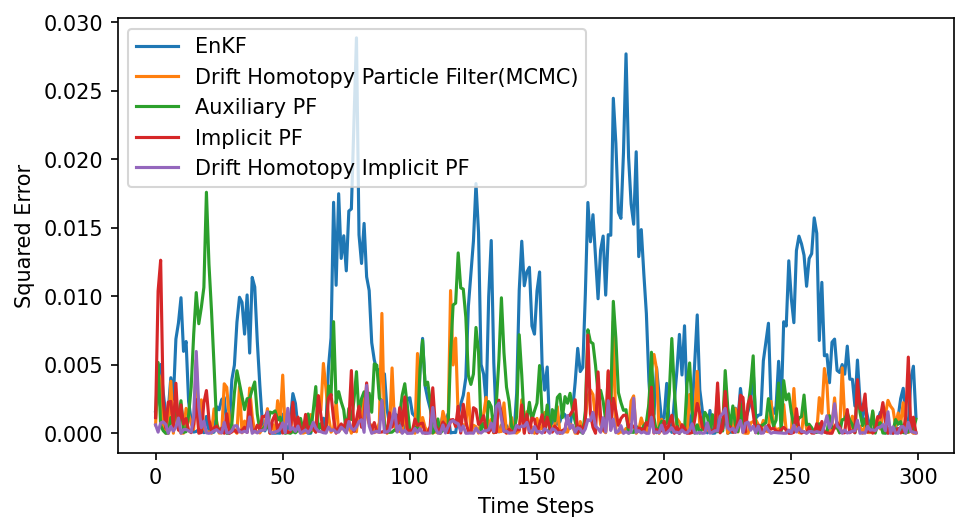} }\end{center}
\caption{Double well potential case 1: $\alpha = 1$, $\sigma = 1.5$, $R = 1.5$ }\label{Ex1_case1}
\end{figure}
In Figure \ref{Ex1_case1} (a), we present the state estimation for the target state obtained by different nonlinear filtering methods, where the black curve is the true simulated state and the colored curves are estimates. From this subplot, we can see that all the methods captured the main trend of the state while the EnKF failed to capture detailed behaviors of the target. To show more details of estimation accuracy, in Figure \ref{Ex1_case1} (b) we plot errors of each method in estimating the state $X$. From this subplot, we can see that the IPF, the DHPF and the DHIPF can give very accurate estimates for the state, the APF occasionally suffers large errors, and the EnKF has the worst performance.

\begin{table}[h!]\small
\leftmargin=6pc \caption{Example 1. Performance comparison for Case 1} \label{e1:c1}
\begin{center}
\begin{tabular}{|c|c|c|c|c|c|c|c|}
\hline   & APF& EnKF & IPF & DHPF & \textbf{DHIPF}   \\
\hline   CPU Time &  $9.703$ & {\color{blue}$0.365625$} & $0.0938$ & $48.563$ & $\bf{0.312}$ \\
\hline  MSE & $2.03 E-3$ & {\color{blue}$5.22 E-3$} & $1.10 E-3$ & $5.92 E-4$ & $\bf{4.60 E-4}$ \\
\hline
\end{tabular}\end{center}
\end{table}
In order to further demonstrate the performance of all the methods, we present the accumulated mean square error (MSE) of each method along with its CPU time in solving this double well potential tracking problem in Table \ref{e1:c1}. From the comparison table, we can see that the DHPF (with MCMC sampling) and the DHIPF have the lowest estimation errors. However, the DHIPF spends much less CPU time compared with DHPF due to the efficient implicit sampling procedure.

\subsubsection*{Case 2.}
\begin{figure}[h!]
\begin{center}
\subfloat[Tracking performance. ]{\includegraphics[scale = 0.55]{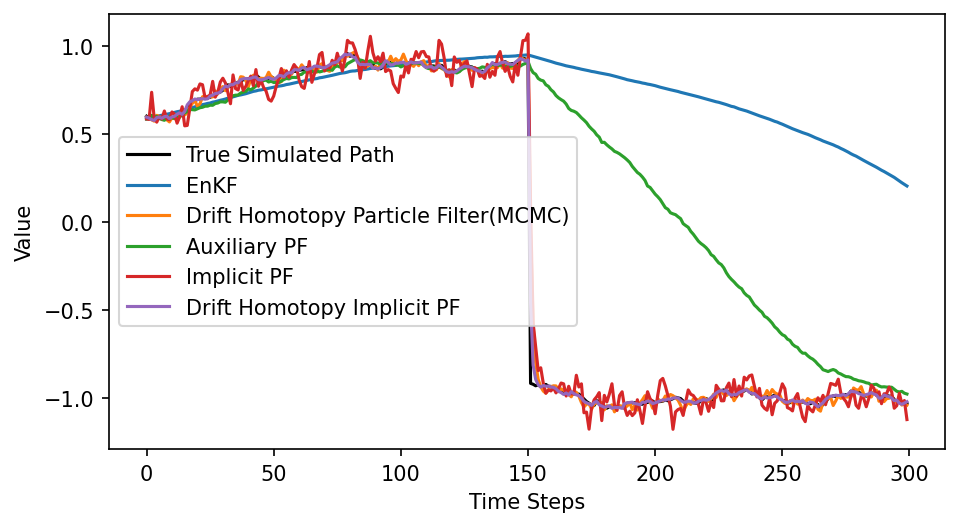} }  \\
\subfloat[Estimation errors.]{\includegraphics[scale = 0.55]{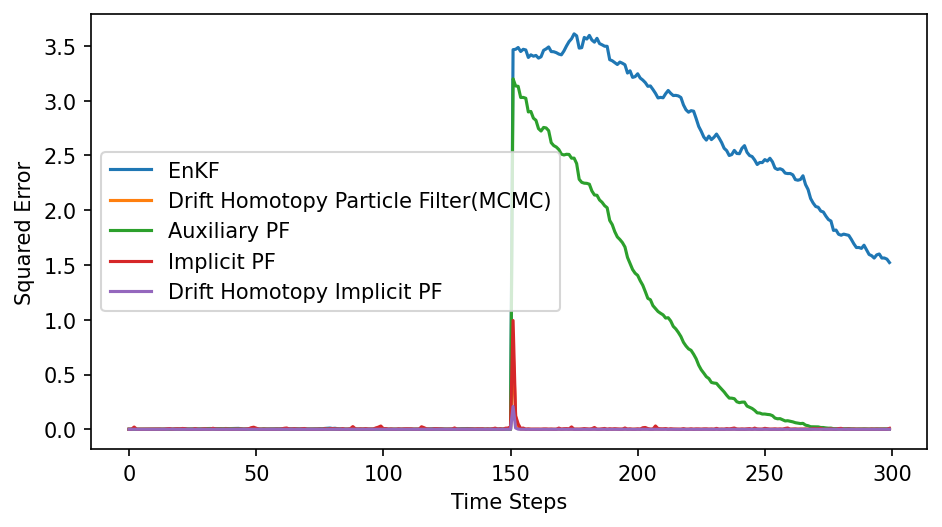} }\end{center}
\caption{Double well potential case 2: $\alpha = 1$, $\sigma = 1$, $R = 1$ with state switch }\label{Ex1_case2}
\end{figure}

In this case, we let $\alpha = 1$, $\sigma = 1$, and $R = 1$, and the initial state $X$ is set to be $X_0 = 0.6$. Different from the first case, we observe that the real state switched from the stable position $X = 1$ to the stable position $X = - 1$ after approximately $150$ tracking steps. In physics, this kind of switch may be caused by some unexpected external force or some extreme diffusion activities.  We use this experiment to demonstrate the ``data first'' advantage of drift homotopy procedure, and we present the tracking performance and the estimation error of each method in Figure \ref{Ex1_case2} (a) and (b), respectively.
From this figure, we can see that the EnKF and APF could not follow effectively the rapid change in the target state. On the other hand, the IPF, the DHPF and the DHIPF captured the switch of the real state effectively while the IPF has slightly higher error at the time of switch.
\begin{table}[h!]\small
\leftmargin=6pc \caption{Example 1. Performance comparison for Case 2} \label{e1:c2}
\begin{center}
\begin{tabular}{|c|c|c|c|c|c|c|c|}
\hline   &   APF& EnKF & IPF & DHPF & \textbf{DHIPF}   \\
\hline   CPU Time & $9.391$ & {\color{blue}$0.578$} & $0.109$ & $43.344$ & $\bf{0.297}$ \\
\hline  MSE & $6.29 E-1$ & {\color{blue}$1.36$} & $5.28 E-3$ & $1.78 E-3$ & $\bf{1.08 E-3}$ \\
\hline
\end{tabular}\end{center}
\end{table}
In Table \ref{e1:c2}, we present the CPU time and the MSE of each method. We can see from this table that due to the unexpected switch of the target state, both the APF and the EnKF have high estimation errors. The DHPF and the DHIPF have very low MSEs in tracking the state, and the DHIPF has much lower computational cost compared with DHPF.

\subsubsection*{Case 3.}
In this case, we also consider state switch during the tracking period. This time, we choose parameters $\alpha = 10$, $\sigma = 1$, and $R = 2$ for the nonlinear filtering problem, where the large parameter $\alpha$ indicates that the potential wells are very ``deep''. 
\begin{figure}[h!]
\begin{center}
\subfloat[Tracking performance. ]{\includegraphics[scale = 0.55]{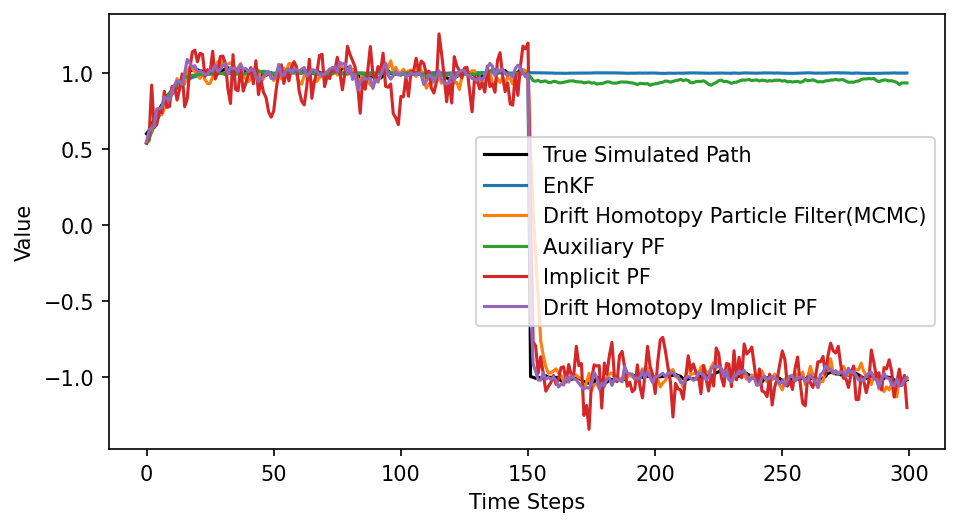} }  \\
\subfloat[Estimation errors.]{\includegraphics[scale = 0.55]{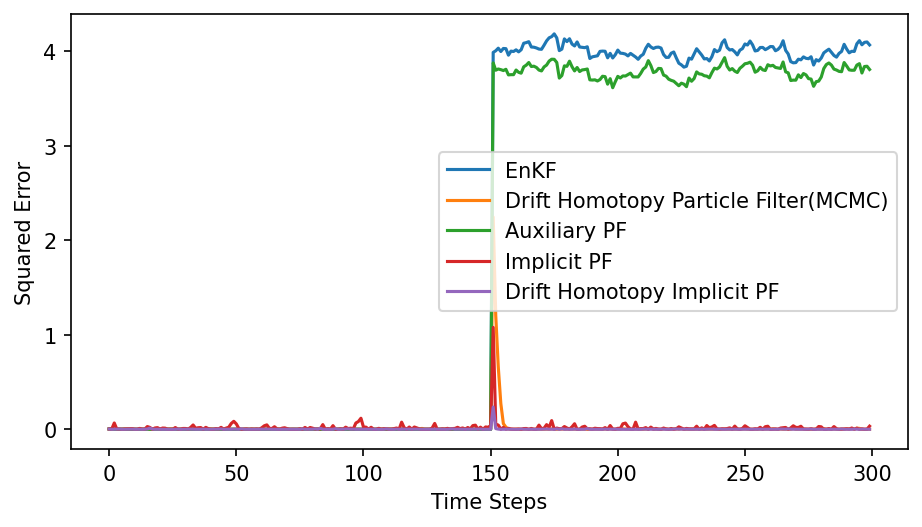} }\end{center}
\caption{Double well potential case 3: $\alpha = 10$, $\sigma = 1$, $R = 2$ with state switch}\label{Ex1_case3}
\end{figure}
In Figure \ref{Ex1_case3}, we present the tracking performance and the estimation error of each method. We can see from this figure that both the EnKF and the APF completely failed to capture the state switch. On the other hand, IPF, DHPF and DHIPF captured the state switch well.  
\begin{table}[h!]\small
\leftmargin=6pc \caption{Example 1. Performance comparison for Case 3} \label{e1:c3}
\begin{center}
\begin{tabular}{|c|c|c|c|c|c|c|c|}
\hline   &   APF& EnKF & IPF & DHPF & \textbf{DHIPF}   \\
\hline   CPU Time & $9.563$ & {\color{blue}$0.453$} & $0.156$ & $50.188$ & $\bf{0.422}$ \\
\hline  MSE & $1.80$ & {\color{blue}$1.99$} & $5.02 E-3$ & $1.35 E-2$ & $\bf{1.59 E-3}$ \\
\hline
\end{tabular}\end{center}
\end{table}
In Table \ref{e1:c3}, we present the CPU time and the MSE of each method. We can see that DHIPF has the lowest accumulative error, and its computational cost is comparable to IPF.

\vspace{0.5em}
From the above numerical experiments, we can see that our DHIPF method outperforms EnKF in accuracy,  and it steadily outperforms APF and DHPF  in both efficiency and accuracy. In comparison with IPF, DHIPF typically has higher accuracy. On the other hand, the computational cost for DHIPF is slightly higher than IPF since DHIPF requires several implicit sampling procedures.
To give more details of advantageous performance of DHIPF compared with IPF, in the next numerical example we focus on the comparison between DHIPF and IPF.

\subsection{Lorenz attractor}

In this example, we solve a Lorenz attractor problem, which has wide applications in weather forecasting and climate prediction. The Lorenz dynamics that we consider is the Lorenz 63 model, which is described by
\begin{equation}\label{Lorenz63}
Lz(\bf{x}) = 
\begin{pmatrix}
a_1(y - x) \\
a_2 x - y  - xz \\
xy  - a_3 z
\end{pmatrix},
\end{equation}
where ${\bf x} = (x, y, z)$ is a three-dimensional vector, $a_1$ is the Prandtl number, $a_2$ is a normalized Rayleigh number and $a_3$ is a non-dimensional wavenumber. The nonlinear filtering problem corresponding to the Lorenz 63 model \eqref{Lorenz63} is given by
\begin{equation}\label{Ex2:NLF}
\begin{aligned}
d X_{t} =& d Lz(t) d t + \sigma dW_t,  \\
d Y_{t} =& X_{t}dt + R d V_t, 
\end{aligned}
\end{equation}
where $W$ and $V$ are two standard Brownian motions, $\sigma$ is the diffusion coefficient, which decides the size of noises that perturb the state process, and $R$ is the coefficient for observational noises. 
In our numerical experiments, we choose $a_1 = 10$, $a_2 = 28$ and $a_3 = 8/3$, which will result chaotic behavior of the state process, and we let $\sigma = {\bf I}_{3\times 3}$, $R = {\bf I}_{3\times 3}$.

\begin{figure}[h!]
\begin{center}
\includegraphics[scale = 0.765]{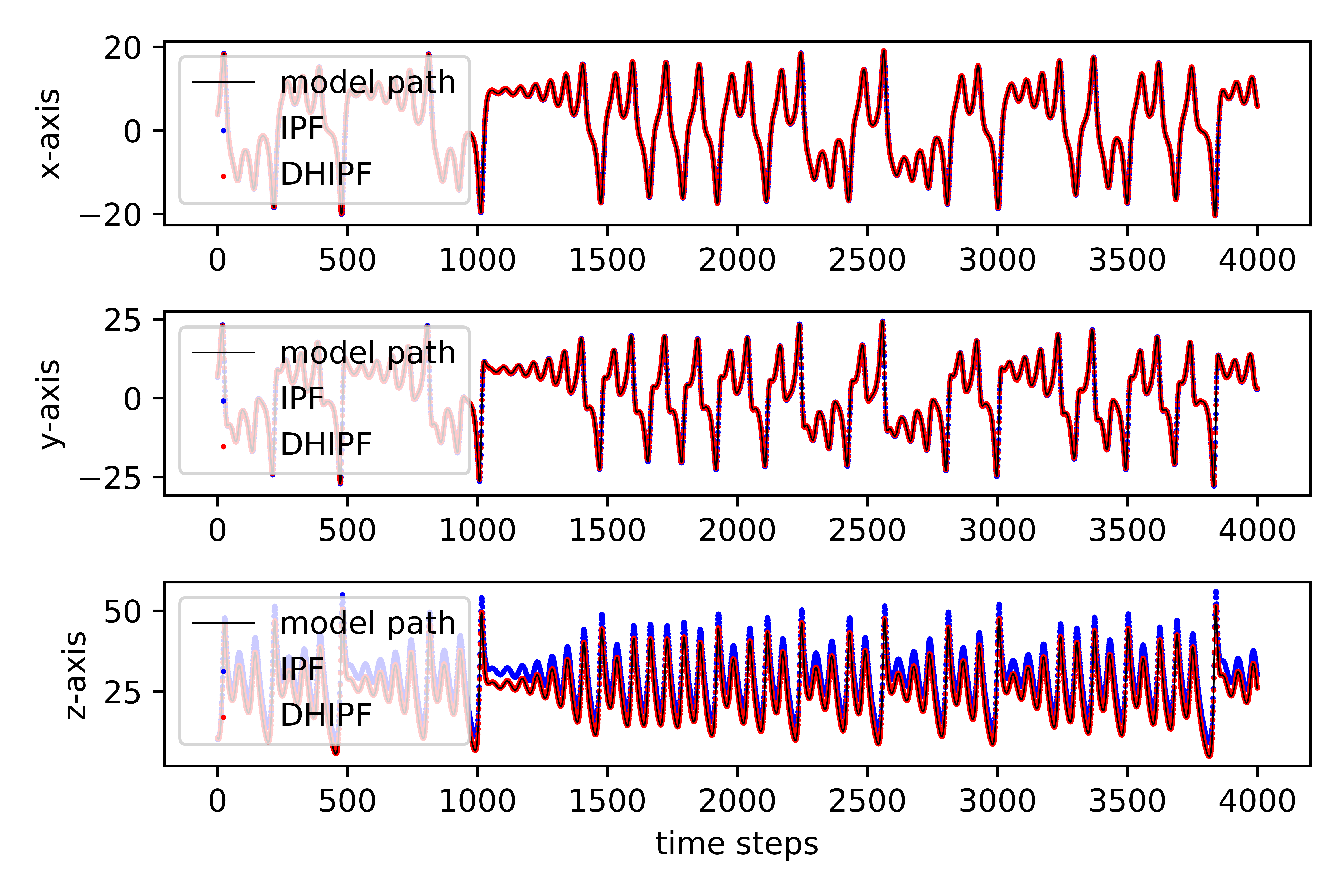} 
\end{center}\vspace{-1em}
\caption{Tracking performance for $4000$ steps.}\label{Ex2_path}
\end{figure}
In the first numerical experiment, we estimate the state $X$ over the time interval $[0, 40]$ with step-size $0.01$, i.e. $4000$ steps, and we use $10$ particles to implement both DHIPF and IPF.  In Figure \ref{Ex2_path}, we present the tracking performance of DHIPF and IPF with respect to each dimension. The black curve is the true target state, the blue dotted curve is the estimate obtained by IPF, and the red dotted curve is the estimate obtained by DHIPF. From this figure, we can see that generally DHIPF and IPF provide good estimates for the target state. Specifically, they both give accurate estimates in $x$ and $y$ directions, and DHIPF is consistently more accurate than IPF in the $z$ direction.  
\begin{figure}[h!]
\begin{center}
\includegraphics[scale = 0.6]{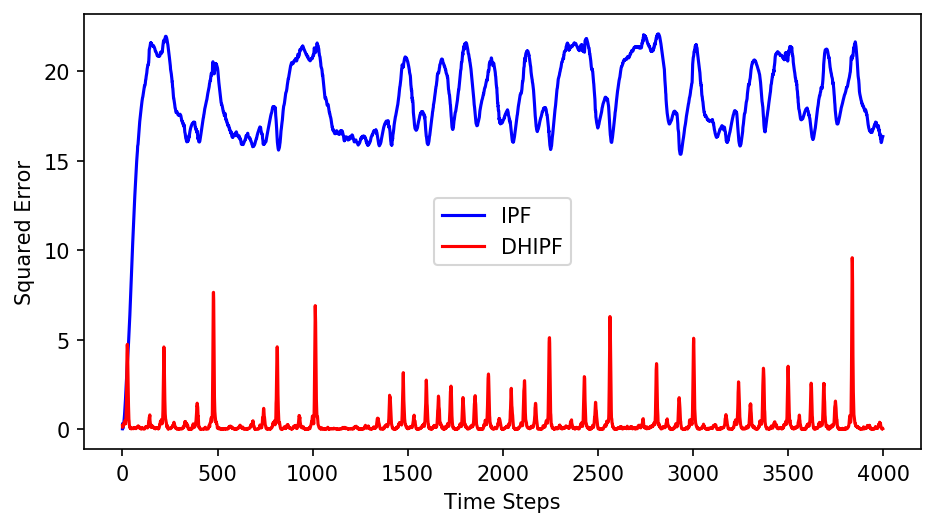} 
\end{center}\vspace{-1em}
\caption{Tracking errors for $4000$ steps. }\label{Ex2_path_error}
\end{figure}
To take a closer look at the accuracy between DHIPF and IPF, we plot squared errors combining all three directions in Figure \ref{Ex2_path_error}. From this figure, wes see that DHIPF is more accurate than IPF  over the entire tracking period.

In our second numerical experiment in this Lorenz attractor example, we present the tracking performance of DHIPF and IPF when the true state of the Lorenz dynamics moves rapidly between time steps $300$ and $400$. Such rapid motion are typically caused by the chaotic nature of the Lorenz model, which is often observed when predicting weather in real time.
To implement DHIPF and IPF, we use $30$ particles to adjust the possible fast state change and we estimate the target state for $600$ time steps.

\begin{figure}[h!]
\begin{center}
\includegraphics[scale = 0.7]{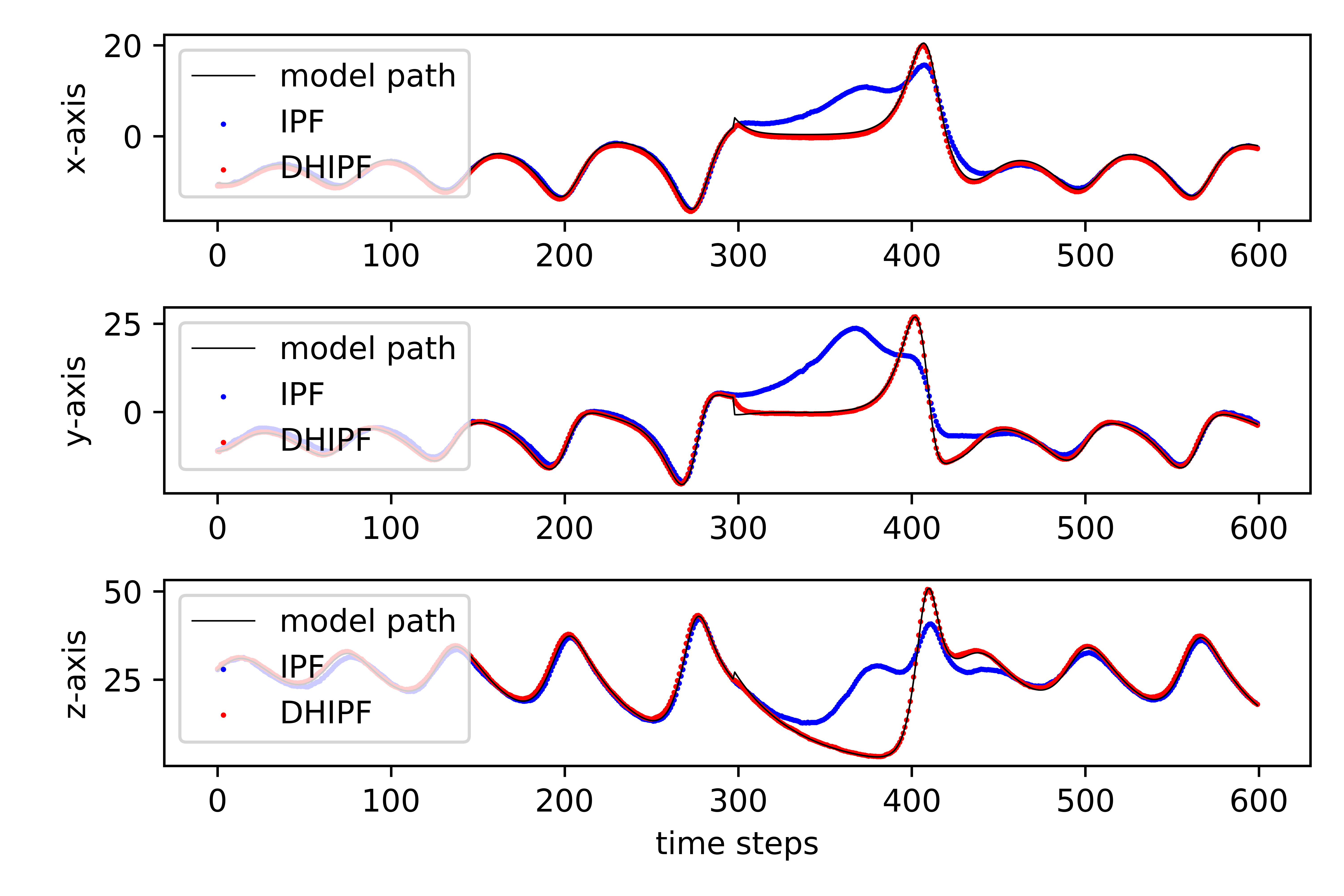} 
\end{center}\vspace{-1em}
\caption{Tracking performance with rapid change in the state. }\label{Ex2_switch}
\end{figure}
In Figure \ref{Ex2_switch}, we show the estimation performance of DHIPF and IPF in each dimension.
We can see from the figure that DHIPF accurately captured the true state of the Lorenz dynamics even during the rapid motion period. On the other hand, IPF could only follow the trend of the target motions and it took IPF over $100$ steps to recover good estimates. 
\begin{figure}[h!]
\begin{center}
\includegraphics[scale = 0.6]{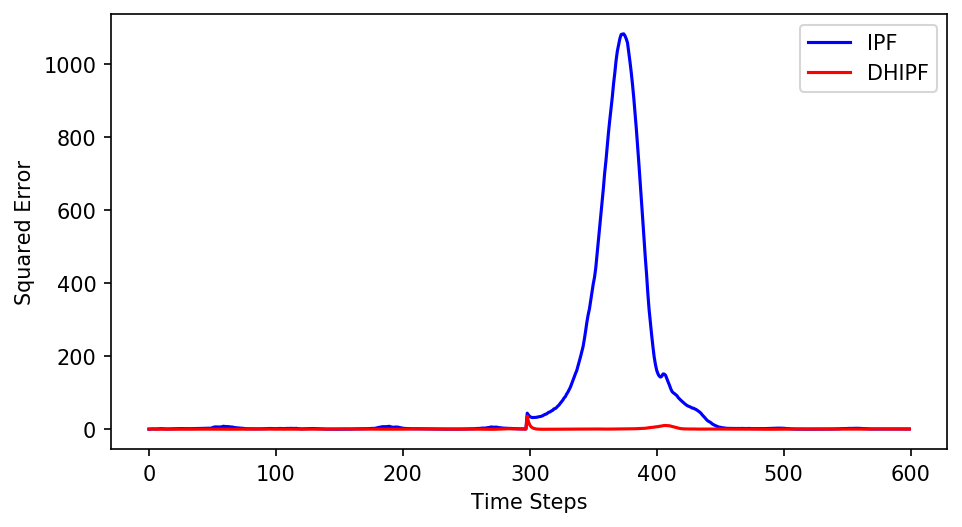} 
\end{center}\vspace{-1em}
\caption{Tracking errors with rapid change in the state. }\label{Ex2_switch_error}
\end{figure}
In Figure \ref{Ex2_switch_error}, we plot the squared errors combining three dimensions over the tracking time, and we can see clearly that IPF suffers large errors from time step $300$ to approximately time step $450$ while DHIPF only has a small spike in estimation errors to adjust the chaotic behavior of the model. 

\vspace{0.4em}

To further demonstrate the advantageous performance of DHIPF over IPF, we assume that there are gaps between model simulations and observations. This could reflect the situation that the data are collected occasionally, and such a situation occurs frequently in practice. In this experiment, we use $50$ particles for both DHIPF and IPF, and we track the target state with $1000$ simulation steps. 
\begin{figure}[h!]
\begin{center}
\includegraphics[scale = 0.6]{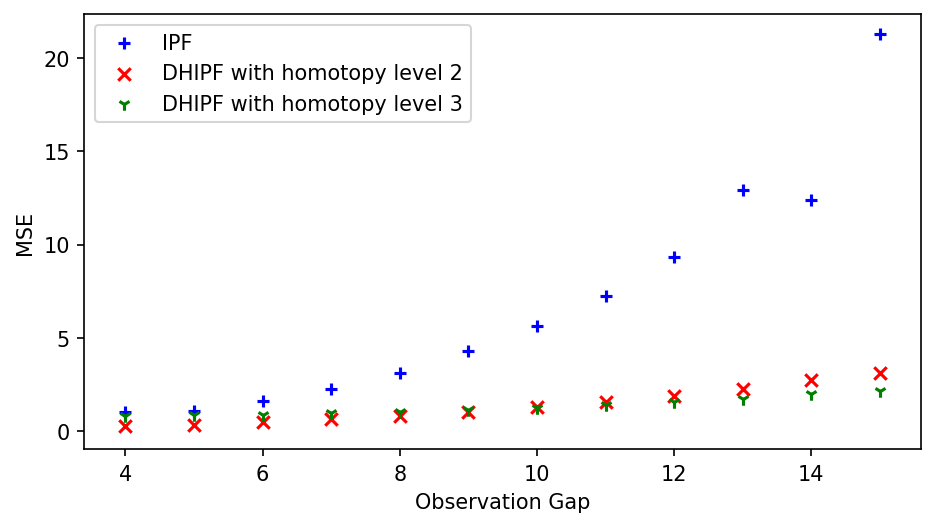} 
\end{center}\vspace{-1em}
\caption{Mean square errors with respect to observation gaps. }\label{Ex2_obs_gap}
\end{figure}
In Figure \ref{Ex2_obs_gap}, we solve the Lorenz attractor problem repeatedly over $20$ times and plot the MSEs among all $20$ repeated tests of each method with respect to observation gaps, where the blue markers are MSEs of IPF, the red markers are MSEs of DHIPF with $2$ drift homotopy levels, i.e. $L=2$, and the green markers are MSEs of DHIPF with $3$ drift homotopy levels, i.e. $L=3$. From this figure, we see that the errors of IPF increase as observations become sparser. On the other hand, although DHIPF has larger errors when observations are sparser, it is not as sensitive as IPF with respect to observation gaps, and DHIPF always has accurate estimates for the state.The reason why DHIPF has better performance in this ``observation gap experiment'' is that the drift homotopy procedure allows us to process the observational data first before we incorporate the dynamical model. As a result, in the case that data are hard to collect, which means each set of observational data is ``more valuable'', DHIPF can utilize the observational data more effectively and therefore obtain better results. 
Moreover, we see from this figure that for smaller observation gaps, DHIPF with $2$ homotopy levels has similar performance to $3$ homotopy levels. When the observation gap is getting larger, more homotopy levels bring more accurate results. This also supports the utility of the homotopy procedure.

\section{Acknowledgement}

This work is partially supported by the Scientific Discovery through Advanced Computing (SciDAC) program funded by U.S. Department of Energy, Office of Science, Advanced Scientific Computing Research through FASTMath Institute and CompFUSE project.  The second author also acknowledges support by U.S. National Science Foundation under Contract DMS-1720222.

\bibliographystyle{plain}


\end{document}